\documentclass{article}

\input{tcilatex}
\begin{document}

\textbf{Transformations of Jesus Guillera's formulas for 1/}$\mathbf{\pi }^{%
\mathbf{2}}$.%
\[
\text{Gert Almkvist} 
\]

\textbf{Introduction.}

Jesus Guillera has found nine formulas of Ramanujan type for \ $1/\pi ^{2}$
\ (see [5] and[6]). The first three are proved using the WZ-method. Then
there exist a lot of formulas obtained by two methods, squaring a formula
for \ $1/\pi $ \ or transforming one of Guillera's formulas. Both methods
were invented by Wadim Zudilin ([7] and[8]). I gave some more examples of
the squaring method in [3]. In [4] Baruah and Berndt obtain a lot of
formulas for \ $1/\pi ^{2}$ by using Eisenstein series, but they are all
examples of squaring formulas for \ $1/\pi $ . Here we give some examples of
very simple transformations. Finally we record some congruences \ $\func{mod}
$ $p^{5}$ not found in Zudilin [9].

\textbf{Transformations.}

Consider the following transformation of sequences%
\[
b_{n}=\sum_{j=0}^{n}(-1)^{j}\binom{n}{j}a_{j} 
\]%
It is an involution since we have the inverse%
\[
a_{n}=\sum_{j=0}^{n}(-1)^{j}\binom{n}{j}b_{j} 
\]%
The corresponding formulas for the generating functions are%
\[
\sum_{n=0}^{\infty }b_{n}x^{n}=\frac{1}{1-x}\sum_{n=0}^{\infty }a_{n}(\frac{%
-x}{1-x})^{n} 
\]%
and%
\[
\sum_{n=0}^{\infty }a_{n}x^{n}=\frac{1}{1-x}\sum_{n=0}^{\infty }b_{n}(\frac{%
-x}{1-x})^{n} 
\]

\textbf{Remark.}

This transformation was used in [1] to find 14 new Calabi-Yau differential
equations, denoted \ $\widehat{1},\widehat{2},...,\widehat{14}$ \ in the
table [2] . Indeed let e.g. ( $(c)_{n}=c(c+1)...(c+n-1)$ \ is the Pochhammer
symbol)%
\[
a_{n}=2^{10n}\frac{(1/2)_{n}^{5}}{n!^{5}} 
\]%
and 
\[
b_{n}=2^{10n}\sum_{j=0}^{n}(-1)^{j}\binom{n}{j}\frac{(1/2)_{j}^{5}}{j!^{5}} 
\]%
Then%
\[
w_{0}=\sum_{n=0}^{\infty }b_{n}x^{n} 
\]%
satisfies a fifth order differential equation whose fourth order pullback is
the Calabi-Yau equation $\widehat{3}$ .

Assume now that we have one of Guillera's formulas for $\ 1/\pi ^{2}$ (see
[5] )%
\[
\sum_{n=0}^{\infty }a_{n}(an^{2}+bn+c)x_{0}^{n}=\frac{1}{\pi ^{2}} 
\]%
Let \ $\theta =x\frac{d}{dx}$ \ and \ 
\[
f(x)=\sum_{n=0}^{\infty }a_{n}x^{n} 
\]%
Then the formula for \ $1/\pi ^{2}$ \ is obtained by substituting \ $x=x_{0}$
\ in%
\[
(a\theta ^{2}+b\theta +c)f(x) 
\]%
This is equal to%
\[
(a\theta ^{2}+b\theta +c)\left\{ \frac{1}{1-x}\sum_{n=0}^{\infty }b_{n}(%
\frac{-x}{1-x})^{n}\right\} 
\]%
\[
=\sum_{n=0}^{\infty }(-1)^{n}b_{n}\left\{ (a\theta ^{2}+b\theta +c)\frac{%
x^{n}}{(1-x)^{n+1}}\right\} 
\]%
Now we have%
\[
\theta \left\{ \frac{x^{n}}{(1-x)^{n+1}}\right\} =\frac{n+x}{(1-x)^{2}}(%
\frac{x}{1-x})^{n} 
\]%
and%
\[
\theta ^{2}\left\{ \frac{x^{n}}{(1-x)^{n+1}}\right\} =\frac{n^{2}+3xn+x+x^{2}%
}{(1-x)^{3}}(\frac{x}{1-x})^{n} 
\]%
It follows%
\[
\sum_{n=0}^{\infty }b_{n}\left\{ a\frac{n^{2}+3x_{0}n+x_{0}+x_{0}^{2}}{%
(1-x_{0})^{3}}+b\frac{n+x_{0}}{(1-x_{0})^{2}}+\frac{c}{1-x_{0}}\right\}
w_{0}^{n}=\frac{1}{\pi ^{2}} 
\]%
where%
\[
w_{0}=\frac{-x_{0}}{1-x_{0}} 
\]

We use Guilleras numbering of his formulas in [5]. The corresponding
Calabi-Yau equations use the notation of [2].

\textbf{Identity 1.}%
\[
a_{n}=\frac{(1/2)_{n}^{5}}{n!^{5}} 
\]%
\[
\sum_{n=0}^{\infty }(-1)^{n}a_{n}(20n^{2}+8n+1)\frac{1}{2^{2n}}=\frac{8}{\pi
^{2}} 
\]

$\widehat{\mathbf{3}}$%
\[
b_{n}=\sum_{j=0}^{n}(-1)^{j}\binom{n}{j}\frac{(1/2)_{j}^{5}}{j!^{5}} 
\]%
\[
\sum_{n=0}^{\infty }b_{n}(64n^{2}-16n-15)\frac{1}{5^{n}}=\frac{50}{\pi ^{2}} 
\]%
\[
\]

\textbf{Identity 2.}

Here \ $a_{n}$ \ and \ $b_{n}$ \ are the same as in Identity 1.%
\[
\sum_{n=0}^{\infty }(-1)^{n}a_{n}(820n^{2}+180n+13)\frac{1}{2^{10n}}=\frac{%
128}{\pi ^{2}} 
\]%
\[
\sum_{n=0}^{\infty }b_{n}(4194304n^{2}+909312n+61633)\frac{1}{1025^{n}}=%
\frac{5^{5}41^{2}}{8\pi ^{2}} 
\]%
\[
\]

\textbf{Identity 3.}%
\[
a_{n}=\frac{(1/2)_{n}^{3}(1/4)_{n}(3/4)_{n}}{n!^{5}} 
\]%
\[
\sum_{n=0}^{\infty }(-1)^{n}b_{n}(1977326743n^{2}+315416969n+16389226)\frac{1%
}{2400^{n}}=\frac{2^{11}3^{2}5^{5}\sqrt{7}}{\pi ^{2}} 
\]

$\widehat{\mathbf{6}}$%
\[
b_{n}=\sum_{j=0}^{n}(-1)^{j}\binom{n}{j}\frac{(1/2)_{j}^{3}(1/4)_{j}(3/4)_{j}%
}{j!^{5}} 
\]%
\[
\sum_{n=0}^{\infty }(-1)^{n}b_{n}(2048n^{2}+928n+215)\frac{1}{15^{n}}=\frac{%
450}{\pi ^{2}} 
\]%
\[
\]

\textbf{Identity 4.}%
\[
a_{n}=\frac{(1/2)_{n}(1/3)_{n}(2/3)_{n}(1/4)_{n}(3/4)_{n}}{n!^{5}} 
\]%
\[
\sum_{n=0}^{\infty }(-1)^{n}a_{n}(252n^{2}+63n+5)\frac{1}{48^{n}}=\frac{48}{%
\pi ^{2}} 
\]

$\widehat{\mathbf{11}}$%
\[
b_{n}=\sum_{j=0}^{n}(-1)^{j}\binom{n}{j}\frac{%
(1/2)_{j}(1/3)_{j}(2/3)_{j}(1/4)_{j}(3/4)_{j}}{j!^{5}} 
\]%
\[
\sum_{n=0}^{\infty }b_{n}(41472n^{2}+7992n-209)\frac{1}{7^{2n}}=\frac{7^{5}}{%
2\pi ^{2}} 
\]%
\[
\]

\textbf{Identity 5.}%
\[
a_{n}=\frac{(1/2)_{n}(1/4)_{n}(3/4)_{n}(1/6)_{n}(5/6)_{n}}{n!^{5}} 
\]%
\[
\sum_{n=0}^{\infty }(-1)^{n}a_{n}(1640n^{2}+278n+15)\frac{1}{2^{10n}}=\frac{%
256\sqrt{3}}{3\pi ^{2}} 
\]

$\widehat{\mathbf{12}}$%
\[
b_{n}=\sum_{j=0}^{n}(-1)^{j}\binom{n}{j}\frac{%
(1/2)_{j}(1/4)_{j}(3/4)_{j}(1/6)_{j}(5/6)_{j}}{j!^{5}} 
\]%
\[
\sum_{n=0}^{\infty }b_{n}(25165824n^{2}+4196352n+201903)\frac{1}{1025^{n}}=%
\frac{5^{5}41^{2}\sqrt{3}}{4\pi ^{2}} 
\]%
\[
\]

\textbf{Identity 6.}%
\[
a_{n}=\frac{(1/2)_{n}(1/8)_{n}(3/8)_{n}(5/8)_{n}(7/8)_{n}}{n!^{5}} 
\]%
\[
\sum_{n=0}^{\infty }a_{n}(1920n^{2}+304n+15)\frac{1}{7^{4n}}=\frac{56\sqrt{7}%
}{\pi ^{2}} 
\]

$\widehat{\mathbf{7}}$%
\[
b_{n}=\sum_{j=0}^{n}(-1)^{j}\binom{n}{j}\frac{%
(1/2)_{j}(1/8)_{j}(3/8)_{j}(5/8)_{j}(7/8)_{j}}{j!^{5}} 
\]%
\[
\sum_{n=0}^{\infty }(-1)^{n}b_{n}(1977326743n^{2}+315416969n+16389226)\frac{1%
}{2400^{n}}=\frac{2^{11}3^{2}5^{5}\sqrt{7}}{\pi ^{2}} 
\]%
\[
\]

\textbf{Identity 7.}%
\[
a_{n}=\frac{(1/2)_{n}(1/3)_{n}(2/3)_{n}(1/6)_{n}(5/6)_{n}}{n!^{5}} 
\]%
\[
\sum_{n=0}^{\infty }(-1)^{n}a_{n}(5418n^{2}+693n+29)\frac{1}{80^{3n}}=\frac{%
128\sqrt{5}}{\pi ^{2}} 
\]

$\widehat{\mathbf{8}}$%
\[
b_{n}=\sum_{j=0}^{n}(-1)^{j}\binom{n}{j}\frac{%
(1/2)_{j}(1/3)_{j}(2/3)_{j}(1/6)_{j}(5/6)_{j}}{j!^{5}} 
\]%
\[
\sum_{n=0}^{\infty }b_{n}(262144000000n^{2}+33528576000n+1402561253)\frac{1}{%
512001^{n}}=\frac{3^{13}7^{5}43^{2}\sqrt{5}}{8000\pi ^{2}} 
\]%
\[
\]

\textbf{Identity 8.}%
\[
a_{n}=\binom{2n}{n}^{2}\sum_{k=0}^{n}\binom{2k}{k}^{2}\binom{2n-2k}{n-k}^{2} 
\]%
\[
\sum_{n=0}^{\infty }a_{n}(36n^{2}+12n+1)\frac{1}{2^{10n}}=\frac{32}{\pi ^{2}}
\]%
\[
b_{n}=\sum_{k=0}^{n}\sum_{j=0}^{n}(-1)^{k}\binom{n}{k}\binom{2k}{k}^{2}%
\binom{2j}{j}^{2}\binom{2k-2j}{k-j}^{2} 
\]%
\[
\sum_{n=0}^{\infty }(-1)^{n}b_{n}(4194304n^{2}+1409024n+121745)\frac{1}{%
1023^{n}}=\frac{3\cdot 11^{3}31^{3}}{32\pi ^{2}} 
\]%
\[
\]

\textbf{Identity 9. }

In [6] Guillera has found a new formula for \ $1/\pi ^{2}$ \ which is very
different from all the others. Let%
\[
a_{n}=\frac{1}{n+\frac{1}{2}}\binom{2n}{n}^{2}\binom{4n}{2n}^{3} 
\]%
Then%
\[
\sum_{n=0}^{\infty }a_{n}(672n^{3}+472n^{2}+78n+\frac{9}{2})\frac{1}{2^{22n}}%
=\frac{64\sqrt{2}}{\pi ^{2}} 
\]%
Observe that \ $y=\dsum\limits_{n=0}^{\infty }a_{n}x^{n}$ \ does not satisfy
a fifth order differential equation (but one of order six).

Here we need 
\[
\theta ^{3}\left\{ \frac{x^{n}}{(1-x)^{n+1}}\right\} =\frac{%
n^{3}+6xn^{2}+(4x+7x^{2})n+x+4x^{2}+x^{3}}{(1-x)^{4}}(\frac{x}{1-x})^{n}
\]%
We have 
\[
b_{n}=\sum_{j=0}^{n}(-1)^{j}\frac{1}{j+\frac{1}{2}}\binom{n}{j}\binom{2j}{j}%
^{2}\binom{4j}{2j}^{3}
\]%
and get the monster formula%
\[
\sum_{n=0}^{\infty }(-1)^{n}b_{n}U(n)\frac{1}{4194303^{n}}=\frac{%
3^{3}(23\cdot 89\cdot 683)^{4}\sqrt{2}}{2^{15}\pi ^{2}}
\]%
where%
\[
U(n)=33056565380087516495872n^{3}+23218343626230634381312n^{2}
\]%
\[
+3836969069974667657216n+221375102329522137953
\]%
\[
\]

B.Gourevich found the following formula 
\[
\sum_{n=0}^{\infty }\frac{(1/2)_{n}^{7}}{n!^{7}}(168n^{3}+74n^{2}+14n+1)%
\frac{1}{64^{n}}=\frac{32}{\pi ^{3}}
\]%
With%
\[
b_{n}=\sum_{k=0}^{n}(-1)^{k}\binom{n}{k}\frac{(1/2)_{k}^{7}}{k!^{7}}
\]%
we find%
\[
\sum_{n=0}^{\infty }(-1)^{n}b_{n}(2097152n^{3}+1130496n^{2}+347776n+64197)%
\frac{1}{63^{n}}=\frac{3^{7}7^{3}}{2\pi ^{3}}
\]%
\[
\]

In [7] Zudilin used a quadratic transformation to obtain two new formulas
for \ $1/\pi ^{2}$ \ from Guillera's identities 1 and 2. Here 
\[
a_{n}=\binom{2n}{n}\binom{4n}{2n}\sum_{k=0}^{n}16^{n-k}\binom{2k}{k}^{3}%
\binom{2n-2k}{n-k} 
\]%
and the first formula is

\textbf{Z 1.}%
\[
\sum_{n=0}^{\infty }a_{n}(18n^{2}-10n-3)\frac{1}{80^{2n}}=\frac{10\sqrt{5}}{%
\pi ^{2}} 
\]%
We have%
\[
b_{n}=\sum_{j=0}^{n}\sum_{k=0}^{n}(-1)^{j}16^{j-k}\binom{n}{j}\binom{2j}{j}%
\binom{4j}{2j}\binom{2k}{k}^{3}\binom{2j-2k}{j-k} 
\]%
and%
\[
\sum_{n=0}^{\infty }(-1)^{n}b_{n}(3276800n^{2}-1818624n-545735)\frac{1}{%
6399^{n}}=\frac{3^{10}79^{3}\sqrt{5}}{2^{7}5^{3}\pi ^{2}} 
\]%
The other formula is

\textbf{Z 2.}%
\[
\sum_{n=0}^{\infty }a_{n}(1046529n^{2}+227104n+16032)\frac{1}{1025^{2n}}=%
\frac{5^{4}41\sqrt{41}}{\pi ^{2}} 
\]%
Our transform is%
\[
\sum_{n=0}^{\infty
}(-1)^{n}b_{n}(76354828515625n^{2}+16569725866875n+1169782053458)\frac{1}{%
1050624^{n}}=\frac{2^{33}3^{7}19^{3}\sqrt{41}}{41^{3}\pi ^{2}} 
\]%
\[
\]

\textbf{Remark.}

After writing this I found Riordan: Combinatorial Identities (in Russian)
and realized that there are an infinity of such transformations. Indeed let%
\[
b_{n}=\sum_{j=0}^{\infty }(-1)^{j}\binom{n+p}{j+p}a_{j} 
\]%
then \ 
\[
\sum_{n=0}^{\infty }a_{n}x^{n}=\frac{1}{(1-x)^{1+p}}\sum_{n=0}^{\infty
}b_{n}(\frac{-x}{1-x})^{n} 
\]

Identity 1 is then transformed to%
\[
\sum_{n=0}^{\infty }b_{n}(64n^{2}-16(2p+1)n+4p^{2}-16p-15)\frac{1}{5^{n}}=%
\frac{5^{p+2}}{2^{2p-1}\pi ^{2}} 
\]%
and Identity 2 to%
\[
\sum_{n=0}^{\infty }b_{n}(2^{22}n^{2}-2^{13}(p-111)n+4p^{2}-4988p+61633)%
\frac{1}{1025^{n}}=\frac{5^{2p+5}41^{p+2}}{2^{10p+3}\pi ^{2}} 
\]%
\bigskip

In particular if \ $p=1/2$ \ we get from Identities 1 and 2 with%
\[
b_{n}=\sum_{j=0}^{\infty }(-1)^{j}\binom{n+1/2}{j+1/2}\frac{(1/2)_{j}^{5}}{%
j!^{5}} 
\]%
the formulas%
\[
\sum_{n=0}^{\infty }b_{n}(64n^{2}-32n-22)\frac{1}{5^{n}}=\frac{25\sqrt{5}}{%
\pi ^{2}} 
\]%
\[
\sum_{n=0}^{\infty }b_{n}(2^{20}n^{2}+2^{10}13\cdot 17n+14785)\frac{1}{%
1025^{n}}=\frac{5^{6}41^{2}\sqrt{41}}{2^{10}\pi ^{2}} 
\]%
\[
\]

If \ $p=-1/2$ \ then we get 
\[
b_{n}=\binom{2n}{n}\sum_{k=0}^{n}(-256)^{-k}\binom{n}{k}\binom{2k}{k}^{4} 
\]%
with%
\[
\sum_{n=0}^{\infty }b_{n}(32n^{2}-3)\frac{1}{20^{n}}=\frac{10\sqrt{5}}{\pi
^{2}} 
\]%
and%
\[
\sum_{n=0}^{\infty }b_{n}(1048576n^{2}+228352n+16032)\frac{1}{4^{n}1025^{n}}=%
\frac{5^{4}41\sqrt{41}}{\pi ^{2}} 
\]%
These formulas resemble Zudilin's formulas in [7].

Another interesting transformation is%
\[
b_{n}=\sum_{k=0}^{n}\binom{2k}{k}a_{n-k} 
\]%
with inverse%
\[
a_{n}=\sum_{k=0}^{n}\left\{ \binom{2k}{k}-4\binom{2k-2}{k-1}\right\} b_{n-k} 
\]%
One finds%
\[
\sum_{n=0}^{\infty }a_{n}x^{n}=\sqrt{1-4x}\sum_{n=0}^{\infty }b_{n}x^{n} 
\]%
With 
\[
a_{n}=\binom{2n}{n}^{5} 
\]%
we get%
\[
b_{n}=\sum_{k=0}^{n}\binom{2k}{k}\binom{2n-2k}{n-k}^{5} 
\]%
and Identities 1 and 2 transform to%
\[
\sum_{n=0}^{\infty }(-1)^{n}b_{n}(2101250n^{2}+842550n+106497)\frac{1}{%
2^{12n}}=\frac{2^{7}5^{2}41\sqrt{41}}{\pi ^{2}} 
\]%
and%
\[
\sum_{n=0}^{\infty
}(-1)^{n}b_{n}(2817520042025n^{2}+618490757170n+44674554281)\frac{1}{2^{20n}}%
=\frac{2^{14}\sqrt{5}13^{3/2}37^{3/2}109^{3/2}}{\pi ^{2}} 
\]%
\[
\]

Another substitution that leads to a quadratic transformation is the
following%
\[
b_{n}=\sum_{k=0}^{[n/2]}\left\{ \binom{n}{k}-\binom{n}{k-1}\right\} a_{n-2k} 
\]%
with the inverse%
\[
a_{n}=\sum_{k=0}^{[n/2]}(-1)^{k}\binom{n-k}{k}b_{n-2k} 
\]%
The generating functions satisfy%
\[
\sum_{n=0}^{\infty }a_{n}x^{n}=\frac{1}{1+x^{2}}\sum_{n=0}^{\infty }b_{n}(%
\frac{x}{1+x^{2}})^{n} 
\]%
Then with%
\[
b_{n}=\sum_{k=0}^{[n/2]}\left\{ \binom{n}{k}-\binom{n}{k-1}\right\} \binom{%
2n-4k}{n-2k}^{5} 
\]%
Identity 1 transforms to%
\[
\sum_{n=0}^{\infty }(-1)^{n}b_{n}U(n)\left\{ \frac{2^{12}}{97\cdot 257\cdot
673}\right\} ^{n}=\frac{(97\cdot 257\cdot 673)^{3}}{2^{21}\pi ^{2}} 
\]%
where%
\[
U(n)=5629498863124500n^{2}+2251797129330760n+281473399652417 
\]%
Identity 2 transforms to%
\[
\sum_{n=0}^{\infty }(-1)^{n}b_{n}U(n)\left\{ \frac{2^{20}}{257\cdot
4278255361}\right\} ^{n}=\frac{(257\cdot 4278255361)^{3}}{2^{33}\pi ^{2}} 
\]%
where%
\[
U(n)=991319172082192724189512500n^{2}+217606647523420455168904220n 
\]%
\[
+15716035651016544248400757 
\]%
\[
\]

\textbf{Congruences.}

In [9] Zudilin observes that Guilleras formulas 1-7 if summed to \ $p-1$ \
where \ $p$ \ is a prime, usually $\ >3$, also produce congruences \ $\func{%
mod}$ $p^{5}$. He proves one of them, the rest being conjectures. We noted
that also the 9-th formula gives a congruence. 
\[
\sum_{n=0}^{p-1}\binom{2n}{n}^{2}\binom{4n}{2n}^{3}\frac{%
672n^{3}+472n^{2}+78n+9/2}{n+1/2}\frac{1}{2^{22n}}\text{ }\equiv 9\QOVERD( )
{2}{p}p^{2}\text{ }\func{mod}\text{ }p^{5} 
\]

The squares of formulas for \ $1/\pi $ \ rarely produce congruences \ $\func{%
mod}$ $p^{5}.$ There are some exeptions when%
\[
a_{n}=\sum_{k=0}^{n}\frac{(1/2)_{k}^{3}(1/2)_{n-k}^{3}}{k!^{3}(n-k)!^{3}} 
\]%
In [3] we find (it is 5.25 in [4])%
\[
\sum_{n=0}^{p-1}a_{n}(441n^{2}+154n+17)\frac{1}{64^{n}}\equiv 17p^{2}\func{%
mod}\text{ }p^{5} 
\]%
In [4] we find%
\[
\sum_{n=0}^{p-1}a_{n}(3n^{2}+n)\frac{1}{4^{n}}\equiv 0\func{mod}\text{ }p^{6}%
\text{ \ \ \ \ \ \ \ \ \ \ \ \ \ \ \ \ (5.22)} 
\]%
\[
\sum_{n=0}^{p-1}(-1)^{n}a_{n}(16n^{2}+16n+5)\equiv 5p^{2}\func{mod}\text{ }%
p^{5}\text{ \ \ \ \ \ \ \ \ \ (5.26)} 
\]%
\[
\sum_{n=0}^{p-1}(-1)^{n}a_{n}(9n^{2}+5n+1)\frac{1}{8^{n}}\equiv p^{2}\func{%
mod}\text{ }p^{5}\text{ \ \ \ \ \ \ \ \ \ (5.28)} 
\]

\newenvironment{litlist}{
  \begin{list}{}{\setlength{\parsep}{0cm}%
                 \setlength{\itemsep}{0cm}%
                 \setlength{\leftmargin}{\parindent}
                 \setlength{\itemindent}{-\parindent}%
                 \raggedright}}{%
  \end{list}%
}

\newpage
\noindent\textbf{References:}
\begin{litlist}
\item
\textbf{1.} G. Almkvist, W. Zudilin, Differential equations, mirror maps and
zeta values, in: Mirror Symmetry V, N. Yui, S.-T. Yau, and J. D. Lewis (eds), AMS/IP Stud. Adv. Math. 38
(International Press \& Amer. Math. Soc., Providence, RI, 2007) 481--515; arXiv: math/0402386
\item
\textbf{2.} G. Almkvist, C. van Enckevort, D. van Straten, W. Zudilin, Tables of Calabi-Yau equations, arXiv: math/0507430.
\item
\textbf{3.} G. Almkvist, Ramanujan-like formulas for \ $1/\pi ^{2}$ \ \'{a}
la Guillera and Zudilin and Calabi-Yau differential equations, Computer Science J. of Moldova, 17 (2009), 1--21.
\item
\textbf{4.} N. D. Baruah, B. C. Berndt, Ramanujan's Eisenstein series and new
hyper\-geometric-like series for \ $1/\pi ^{2}$, J. Approx. Theory, 2008
\item
\textbf{5.} J. Guillera, My formulas for $1/\pi ^{2}$, available at
\texttt{http://personal.auna.com/guillera}
\item
\textbf{6.} J. Guillera, Expansions related to Ramanujan series and
alike, available at \texttt{http://personal.auna.com/guillera}
\item
\textbf{7.} W. Zudilin, Quadratic transformations and Guillera's formulae,
arXiv CA/0509465
\item
\textbf{8.} W. Zudilin, More Ramanujan-type formulae for \ $1/\pi ^{2}$ ,
Russian Math. Surveys 62 (2007), 634--636.
\item
\textbf{9.} W. Zudilin, Ramanujan-type congruences, arXiv, math.NT0 805.2788%
\end{litlist}

\bigskip\bigskip

\noindent
Institute of Algebraic Meditation\\[0.5mm]
Fogdar\"{o}d 208\\[0.5mm]
SE 243 33 H\"{O}\"{O}R\\[0.5mm]
Sweden\\[0.5mm]
\texttt{gert.almkvist@yahoo.se}

\end{document}